# Data-driven Method for 3D Axis-symmetric Object Reconstruction from Single Cone-beam Projection Data

**Shousheng Luo**

School of Mathematics and Statistics, Data analysis Technology Lab of Henan University, Henan University, Kaifeng, China
sluo@henu.edu.cn

**Ruyue Meng**

School of Mathematics and Statistics, Henan University, Kaifeng, China
991539550@qq.com

**Suhua Wei**

Institute of Applied Physics and Computational Mathematics, Beijing, China
Wei_suhua@iapcm.ac.cn

**Jian-Feng Cai**

Department of Mathematics, Hong Kong University of Science and Technology, Hong Kong
jfcai@ust.hk

**Xue-Cheng Tai**

Department of Mathematics, Hong Kong Baptist University, Hong Kong
xuechengtai@hkbu.edu.hk

**Yang Wang**

Department of Mathematics, Hong Kong University of Science and Technology, Hong Kong
yangwang@ust.hk

**ABSTRACT:** In this paper we consider 3D axis-symmetric (AS) object reconstruction from single cone-beam x-ray projection data. Traditional x-ray CT fails to capture fleeting state of material due to the long time for data acquisition at all angles. Therefore, AS object is devised to investigate the instant deformation of material under pulse change of environment because single projection data is enough to reconstruct its inner structure. Previous reconstruction methods are layer by layer, and ignore the longitudinal tilt of x-ray paths. We propose a regularization method using adaptive tight frame to reconstruct the 3D AS object structure simultaneously. Alternating direction method is adopted to solve the proposed model. More importantly, a numerical algorithm is developed to compute imaging matrix. Experiments on simulation data verify the effectiveness of our method.

## 1. Introduction

X-ray CT (XCT) technique is widely used in various areas, such as medicine, industry detection and material science [8]. Generally speaking, projection data at all angles from 0 to $\pi$ are necessary to reconstruct the attenuation map of object. Therefore, traditional XCT technique fails to capture the instant state of fleeting change of object, for example, the object state under powerful shocks caused by explosion. In this paper, we consider axial-symmetric object reconstruction using single cone-beam projection data, (see [5] for details).

Commonly used methods for the considered problem are the FDK technique [7], i.e. reconstruct the 3D attenuation map slice by slice, which ignores the longitudinal tilt of x-ray paths. For each slice reconstruction, one needs to solve the Abel transform [1, 10]. The reconstruction map by FDK method often suffers from strong artefacts for long objects. Regularization methods are adopted for the problem [2, 5]. In [2], total variation (TV) minimization is adopted to regularize the inversion of Abel transform. In [5], high order total variation regularization is used to reconstruct 3D cylindrical symmetric object from cone-beam projection data slice by slice.

We propose a regularization model based on adaptive tight frame for 3D AS object reconstruction from single cone-beam projection data. Adaptive tight frame is widely used for various image processing problems such as image denoising [3], XCT image reconstruction [11, 12] and PET-MRI [6] due to its high performance in preserving image details and edges. In this paper, we adopt adaptive tight frame regularization to reconstruct 3D AS map simultaneously.

An efficient algorithm is proposed to compute the imaging matrix, i.e. the intersection lengths of all x-ray path with all 3D annular cylinders, for this problem, which makes our method practicable. As far as we know, this is the first work to investigate the 3D AS object reconstruction simultaneously and the computation method for intersection length of a line with annular cylinders. The imaging matrix plays a critical role for the simultaneous reconstruction method. Discretizing the reconstruction cylinder region, which contains the AS object, by annular cylinder element, we can approximate the problem as a linear system, and develop an efficient algorithm to compute the imaging matrix.

Alternating direction method is used to solve the model. There are three variables to update in each loop. We have efficient algorithm for each sub-problem. Firstly, the adaptive frame is learnt from the current map estimate, and then we obtain the representation coefficient by hard thresholding operation, i.e. we have a regularization version of the current estimate map. Lastly, we update the attenuation map by solving a minimization problem with two quadratic terms.

The rest of this paper is organized as follows. In section 2, the proposed model and algorithms are presented. Section 3 is devoted to the computation method for the imaging matrix. Numerical results on simulated data are illustrated in section 4. Conclusions and future work are discussed in section 5.

## 2. The proposed method

In this section, a brief introduction to adaptive tight frame is presented. Then the considered problem and the proposed method are described.

### 2.1 Adaptive wavelet tight frame

A wavelet tight frame is defined by a bank $b = \{b_i\}_{i=1}^{r^2}$ where $b_i \in \mathbb{R}^{r \times r}$ [3]. The analysis operator based on $\{b_i\}_{i=1}^{r^2}$ is defined as follows:

$$W(b): g \in \mathbb{R}^{m \times n} \to v = (b_i * g) \in \mathbb{R}^{r \times m \times n} \tag{1}$$

with $v_i = b_i * g$, where $g$ is an input image. And the synthetic operator is denoted by $W^T(b)$. Hereafter, we omit the option $\{b_i\}_{i=1}^{r^2}$ or $b$ for simplicity, and $\tilde{}$ denote the vector form of a two dimension array by lexicographic order. For a given image $g$, the associating adaptive wavelet tight frame is obtained by minimizing

$$\arg\min_{b,v} \|v - Wg\|_2^2 + \gamma^2 \|v\|_0, W^T(b)W(b) = I, \tag{2}$$

where $I$ is the identity operator in $\mathbb{R}^{m \times n}$ or $\mathbb{R}^N (N = mn)$, and $\gamma > 0$ is a parameter. This problem is solved by alternating direction method as follows:

1. Update $v$: $v^{k+1} = \arg\min_v \|v - W^k g\|_2^2 + \gamma \|v\|_0$.

2. Update filter bank $\{b_i\}_{i=1}^{r^2}$:

$$b^{k+1} = \arg\min_b \|v^k - Wg\|_2^2, W^T W = I. \tag{3}$$

For the first step, $v^{k+1}$ can be obtained by hard thresholding

$$[H_\gamma v](i,j) = \begin{cases} v(i,j) & \text{if } |v(i,j)| > \sqrt{\gamma}, \\ 0 & \text{otherwise}, \end{cases} \tag{4}$$

for $i = 1, 2, \cdots r^2, j = 1, 2, \cdots, mn$.

The solution for the second step is reviewed as follows (see [3] for details). For each $v(i,j)$, $i = 1, 2, \cdots, r^2, j = 1, 2, \cdots, N$, there is a patch of $g$, denoting by $g_j$, such that $v(i,j) = \tilde{g}_j^T \tilde{b}_i$. Let $V = (\tilde{v}_1^T; \tilde{v}_2^T; \cdots; \tilde{v}_{r^2}^T)$. Then we have

$$V = BG, \tag{5}$$

where $B = (\tilde{b}_1^T; \cdots \tilde{b}_{r^2}^T)$, $G = (\tilde{g}_1, \cdots, \tilde{g}_N)$. The condition $W^T W = I$ is guaranteed by $B^T B = \frac{1}{r^2} E$ [3], where $E \in \mathbb{R}^{r^2 \times r^2}$ is the identity matrix. Therefore, the update of $b$ can be reformulated as

$$\{b^{k+1}\} = \arg\min_b \|V^{k+1} - BG\|_2^2, B^T B = \frac{1}{r^2} E, \tag{6}$$

and we have that

$$\tilde{b}_i^{k+1} = \frac{1}{r^2} [U_R U_L]_i, i = 1, 2, \cdots, r^2, \tag{7}$$

where $V^{k+1} G^T = U_L \Sigma U_R$ is the singular value decomposition of $V^{k+1} G^T$, and $[\cdot]_i$ denotes the $i$th column of a matrix.

## 2.2 Problem description

We consider the following imaging system. X-ray source is at $(x_0, 0, 0)$, the detector plane (array) is at $x = x_1$, and the imaging object is enclosed in a three dimensional cylinder domain

$D = \{(x,y,z) \mid \sqrt{x^2+y^2} \leq R_0, |z| \leq Z_0\}$ (see figure 1). Assume the object is axial-symmetric with respect to $z$ axis.

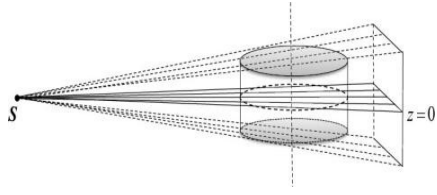

**Figure 1: Imaging system**

Therefore, the attenuation map, denoted by $f$, is an AS function, i.e. there is a two-variable function $u(\cdot,\cdot)$ such that

$$f(x,y,z) = u(\sqrt{x^2+y^2}, z). \tag{8}$$

Discretizing the reconstruction domain $D$ as

$$D_{ij} = \{(x,y,z) \mid r_i \leq \sqrt{x^2+y^2} < r_{i+1}, z_j \leq z < z_{j+1}\}, \tag{9}$$

where $r_i = i\Delta r$, $z_j = -j\Delta z$, $i = 1, 2, \cdots, m$, $j = -n, \cdots, n$ with $m = \left\lceil \frac{R_0}{\delta_r} \right\rceil$, $n = \left\lceil \frac{Z_0}{\delta_z} \right\rceil$, and $\Delta r > 0$, $\Delta z > 0$ are the step sizes.

Suppose the detector array is index by $s = -p, \cdots, 0, 1, p-1$ and $t = -q, \cdots, 0, 1, q-1$. Hereafter, we also denote the radial function $f$ and detector array as a vector by lexicographic order, and $a_{kl}$ (resp. $a_{st,ij}$) denotes intersection length between x-ray path index by $k$ (resp. $(s,t)$) with the l-th (resp. $(i,j)$) annular cylinder elements with $k = 1, 2, \cdots 4pq$, $l = 1, 2, \cdots, 2mn$ ( resp. $s = 1, 2, \cdots p$, $t = 1, 2, \cdots, q, i = 1, 2, \cdots, m, j = 1, 2, \cdots, n)$.

Let $A = (a_{kl})$, and $g = (g_1, g_2, \cdots, g_M)^T$ with $g_k$ being the line integral along the k-th x-ray path determined by the corresponding detector and $M = 4pq$. In practice, $g$ is computed from the recorded x-ray energies before and after placing the object by using Beer's law. According to the aforementioned notations, the algebraic model for the concerned problem can be formulated as

$$Au = g. \tag{10}$$

The computation of $A$ will be discussed in section 3. The reconstruction will suffer strong artifacts because of photon noise and the ill-posedness of $A$.

## 2.3 Data-driven model for AS object reconstruction

In order to suppress the noise effect, we propose a data driven method for AS objection reconstruction as follows

$$\min_{u,b} \|Au - g\|_2^2 + \gamma_1^2 \|W(b)u\|_0, \quad W^T(b)W(b) = I, \tag{11}$$

where $\gamma_1 > 0$ is an user-specified regularization parameter. Our method is different from the one in [12], which needs a reference image to learn the filter bank $b$ previously. In our work, the adaptive tight frame is learnt in the image reconstruction procedure. Therefore, our method does not need any reference image. Introducing an auxiliary variable $v = Wu$, we can relax (11) as

$$\min_{u,b,v} \|Au - g\|_2^2 + \lambda \left[ \|Wu - v\|_2^2 + \gamma_2^2 \|v\|_0 \right], W^T W = I, \tag{12}$$

where $\lambda > 0$ is a penalty parameter, and $\gamma_2^2 = \gamma_1^2 / \lambda$. We use alternating direction method to solve (12), i.e. only one variable among $b, v$ and $u$ is updated with others fixed. The update algorithms for $b$ and $v$ are reviewed in section 2.1. For the update of $u$, we should solve the minimizing problem

$$u^{k+1} = \arg\min_u \frac{1}{2} \|Au - g\|_2^2 + \lambda \|W(b^{k+1})u - v^{k+1}\|_2^2. \tag{13}$$

This problem is equivalent to solve [3]

$$\arg\min_u \frac{1}{2} \|Au - g\|_2^2 + \lambda \|W^T(b^{k+1})v^{k+1} - u\|_2^2. \tag{14}$$

Here we use primal-dual method [4] to solve it, which is presented in Algorithm 1. Details are not given here due to page limit. One can obtain this algorithm using the technique in [4] and [9].

**Algorithm 1.** Algorithm for problem (14)

Let $L \leftarrow \|A\|_2$, $\tau \leftarrow 1/L$, $\sigma \leftarrow 1/L$, $\theta \leftarrow 1$, and set the maximum iteration $N_1$.
(2) Initialize $u^0 = 0, \omega^0 = 0, \bar{u}^0 = u^0$ and $l = 0$.
(3) While $l < N_1$
(4) $\quad \omega^{l+1} = (\omega^l + \sigma(A\bar{u}^l - g))/(1 + \sigma)$;
(5) $\quad u^{l+1} = (2\lambda\sigma W^T(b^{k+1})v^{k+1} + u^l - \tau A^T \omega^{l+1})/(2\lambda\sigma + 1)$;
(6) $\quad \bar{u}^{l+1} = u^{l+1} + \theta(u^{l+1} - u^l)$;
(7) $\quad l = l + 1$;
(8) End (While).
(9) Output $u^{N_1}$ as $u^{k+1}$.

According to the discussions above, we can summarize the algorithm for (12) as Algorithm 2.

## 3. Computation for imaging matrix

According to the assumption, all x-ray paths associating with detector array and x-ray source can be denoted as $L(\gamma, \alpha)$

$$\begin{cases} x = x_0 + t\cos\gamma\cos\alpha, \\ y = t\cos\gamma\sin\alpha, \\ z = t\sin\gamma, \end{cases} \tag{15}$$

where $\gamma$ is the angle between the path and OXY plane, α is the angle between the projection of path on OXY plane and x-axis. Because of the symmetries of the imaging system and measure object, we need to compute the intersection length between the x-ray path with and the annular elements $\gamma, \alpha \geq 0$, and others can be obtained by symmetric property.

---
**Algorithm 2** Algorithm for problem (12)
---
(1) Set maximum iteration number $N_2$, tolerance $\varepsilon$.
(2) Initialize $u^0$ by TV regularization method and $k = 0$.
(3) While $k < N_2$
(4)  Learnt filter bank $b^{k+1}$ by the method in section 2.1;
(5)  Update $v^{k+1}$ by hard thresholding (4);
(6)  Update $u^{k+1}$ by Algorithm 1;
(7)  $k = k + 1$;
(8) End (While)
---

For a given x-ray path $L(\alpha, \gamma)$, the intersections with all annular cylinders are computed in two steps. Firstly, the intersection points of the line and all cylinders with radius $\left\{ (x, y, z) \mid \sqrt{x^2 + y^2} = i\Delta r, z \in \mathbb{R} \right\}$ are

$$\left( x_0 + t_{is} \cos \gamma \cos \alpha, t_{is} \cos \gamma \sin \alpha, t_{is} \sin \gamma \right), \qquad (16)$$

where $i = N(\alpha), N(\alpha) + 1, \cdots, N, s = 0, 1$. The intersection lengths of this ray with the annular cylinders

$$D_{ij} = \left\{ (x, y, z) \mid r_{i-1} < \sqrt{x^2 + y^2} < r_i, z_j < z < z_{j+1} \right\}$$

are computed easily. Firstly, if $i < N(\alpha)$, the intersection length must be zero. We can update the imaging matrix $A$ by Algorithm 3 for any given two intersection points with the annular cylinder. Here we initialize $A$ as zero matrix because it is possible that the intersection length consists of two sections when $\gamma$ is small enough.

---
**Algorithm 3** Update $A$ for given two intersection points between k-th x-ray path and i-th annular cylinder
---
(1) Input two intersection heights $h_1, h_2 (h_1 < h_2)$.
(2) Let $\lfloor p = h_1 / \Delta z \rfloor, q = \lfloor h_2 / \Delta z \rfloor$.
(3) If $p = q$
(4)  $a_{k, p_i} = a_{k, p_i} + |h_2 - h_1| / \sin \gamma$.
(5) Else
(6)  $a_{k, p_i} = a_{k, p_i} + |(p+1)\Delta z - h_1| / \sin \gamma$;
(7)  $a_{k, q_i} = a_{k, q_i} + |h_2 - q\delta z| / \sin \gamma$;
(8)  $a_{k, t_i} = a_{k, t_i} + \delta z / \sin \gamma, p < t < q$;
(9) End (If).
---

In Algorithm 3, $p_i$ (resp. $q_i, t_i$) are the annular index determined by $p$ and $i$. According to Algorithm 3, we can compute the imaging matrix row by row as Algorithm 4. Here $A_k$ denotes the k-row of A in Algorithm 4.

## 4. Numerical experiments

In this section, we will present experiments on simulation data to verify the effectiveness of the proposed method and algorithm for the considered problem, AS object reconstruction from single cone-beam projection. A cylinder phantom with radius 1 and height 2 is designed to test the proposed algorithm (see Figure 2 for the central longitudinal section). For the data acquisition, we set the x-ray source at $(40,0,0)$ and the object center at origin, and the detector array in the plane $x=-50$ with $-2.51 \leq z \leq 2.51, -2.45 \leq y \leq 2.45$.

For the discretization, we set $\Delta z = \Delta r = 0.005$, and the detector size is 0.005×0.005. The imaging matrix is computed by Algorithm 4, and then we generate the projection data by $g = Au + n$, where n is data noise. In this paper, noise obeys Gaussian distribution with mean 0 and variance 0.03. In Algorithm 1, the parameters $\sigma$ and $\tau$ equal to 0.2. The iteration numbers $N_1 = 5000, N_2 = 3$ in Algorithm 1 and Algorithm 2. The proposed method is compared with TV regularization method.

---

**Algorithm 4.** Imaging matrix computation

---

(1) Initializing $A = 0 \in \mathbb{R}^{k,mn}$.
(2) For all x-ray line $L(\gamma_i, \alpha_j)$ numbered by $k$.
(3) Compute the intersection points $t_{i0}$ and $t_{i1}$ of line $L(0,\alpha)$ with the circle with radius $i\Delta r, i = N, N-1, \cdots, N(\alpha)$.
(4) Compute the heights of each intersection points $h_{is} = t_{is} \sin \gamma$
(5) For $i = N(\alpha), \cdots N-1$ and $s = 0,1$;
(6)    Update $A_k$ using algorithm 3 with $h_1 = h_{i0}, h_2 = h_{(i+1)0}$;
(7)    Update $A_k$ using algorithm 3 with $h_1 = h_{i1}, h_2 = h_{(i+1)1}$;
(8) End (For).
(9) Update $A_k$ using algorithm 3 with $h_1 = h_{N(\alpha)0}, h_2 = h_{N(\alpha)1}$.

---

**Table 1: RMSEs of the images in Figure 3**

| λ | TV | ATF |
|---|---|---|
| 0.005 | 0.0323 | 0.0278 |
| 0.01 | 0.0313 | 0.0303 |
| 0.015 | 0.0335 | 0.0325 |

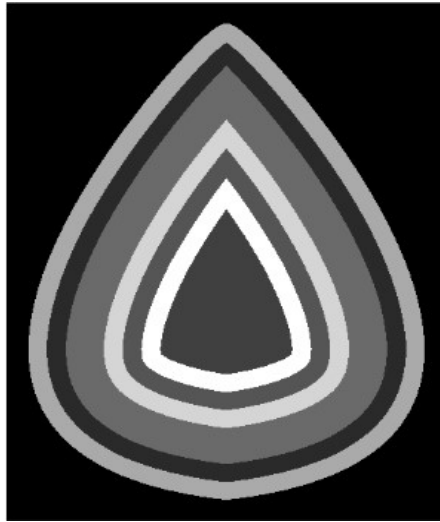

**Figure 2: Central longitudinal section image of the phantom.**

Besides visual quality, we also compare the relative means square error, which is defined as follows:

$$RMSE = \sqrt{\frac{1}{mn}\sum_{i=1}^{m}\sum_{j=1}^{n}(u_{ij}-u_{ij}^{*})^{2}}, \tag{17}$$

where u is the result of the reconstruction, $u^*$ is the original image. The central longitudinal sections of reconstructions with $\lambda = 0.005, 0.01, 0.015$ are illustrated in Figure 3. In order to show the performance of the proposed method, we also illustrate the section images at three different heights in Figure 4. These section images are obtained by rotating the associated profiles of each longitudinal images in Figure 3. The profiles of the reconstructions by the proposed method and TV at $z=0$ are plotted in Figure 5. The RMSEs for different reconstructions are tabulated in Table 1.

It is obviously that the proposed method is superior to the TV regularization method in preserving image edges and suppressing artifacts and in terms of RMSE. As we can see that there are strong artifacts in the reconstructions by TV method, especially when the section is far from the central slice. However, the proposed method suppresses the artifacts effectively. On the other hand, the RMSE values in Table 1 also show the proposed method is better than TV method.

## 5. Conclusion and future work

There are two main contributions in this work. Firstly, this paper presented a data driven simultaneous method for 3D AS object reconstruction from single cone-beam data. Secondly, efficient algorithm for the imaging matrix computation is investigated. In the future, we will study on rapid convergence algorithm for the update of u, e.g. the acceleration of primal-dual algorithm based on the property of the imaging matrix.

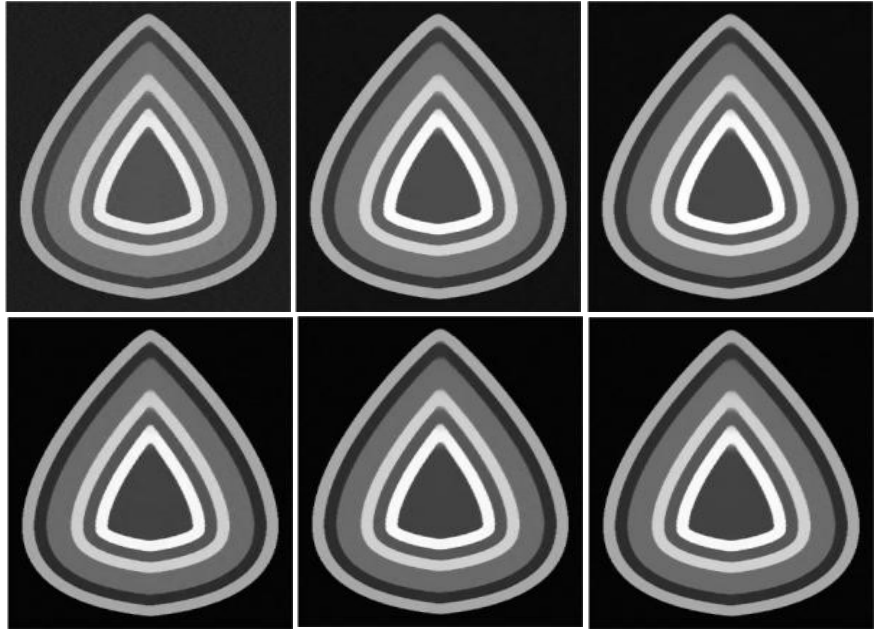

**Figure 3: Central longitudinal section. The reconstructions by TV and the proposed methods are shown in first and second rows, respectively. From left to right λ = 0.005, 0.01 and 0.015 respectively.**

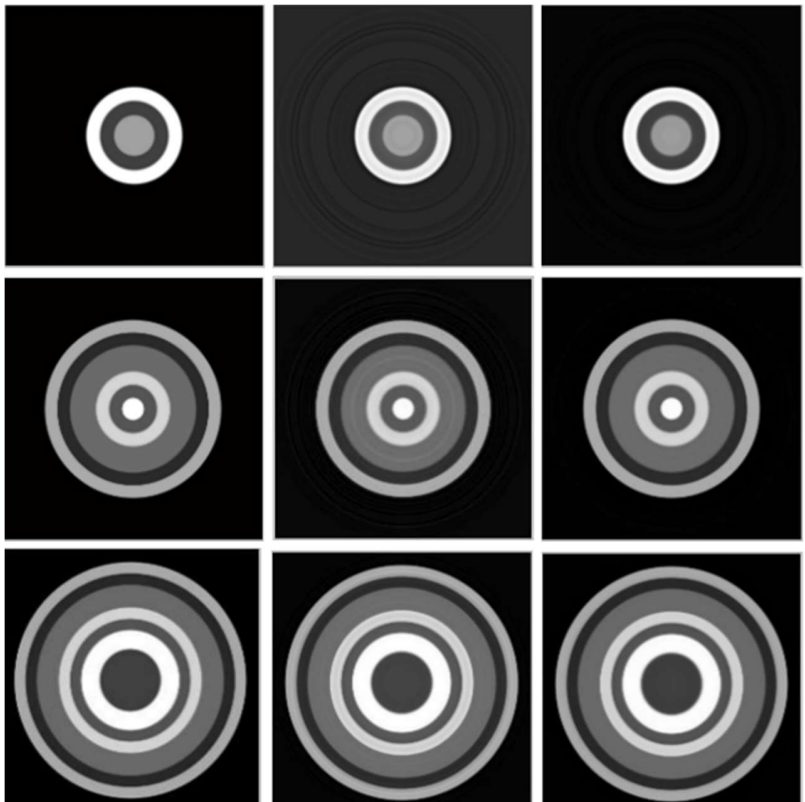

**Figure 4: Images from up to bottom are the original, TV and ATF reconstructions of 100th, 200th and 360th rows. The regularization parameter λ = 0.005, 0.01, 0.015 from left to right.**

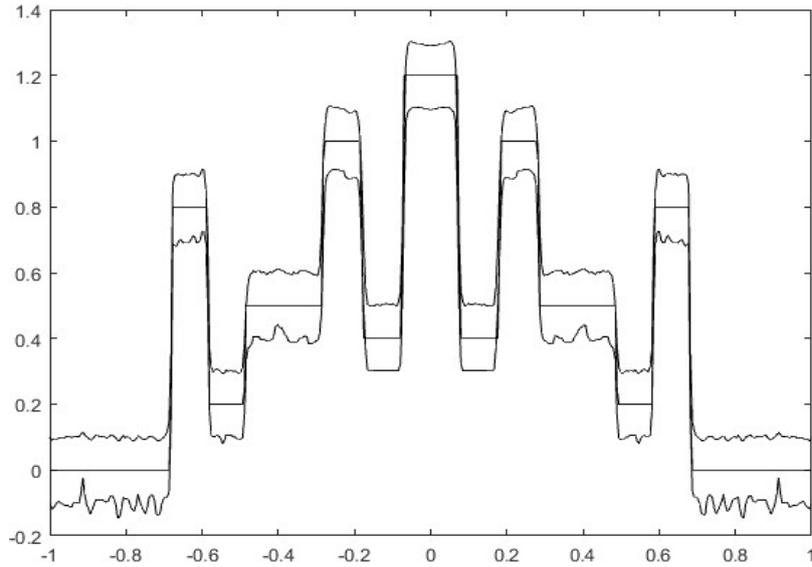

Figure 5: Profile comparison: the reconstructions at z = 0 is plotted. For visual comparison, 0.1 (resp. −0.1) is added to the reconstruction value by the proposed method (resp. TV).

## 6. Acknowledgments

Shousheng Luo was supported by the Programs for Science and Technology Development of He'nan Province (192102310181). Suhua Wei was supported by National Natural Science Foundation of China (11571003). Xue-Cheng Tai was supported by the startup grant at Hong Kong Baptist University, grants RG(R)- RC/17-18/02-MATH and FRG2/17-18/033. Jian-Feng Cai was supported in part by Hong Kong Research Grant Council grant 16306317. Yang Wang was supported in part by the Hong Kong Research Grant Council grants 16306415 and 16308518.